\documentclass[11pt]{amsart}
\usepackage{geometry}                
\usepackage{graphicx}
\usepackage{amssymb}
\usepackage{epstopdf}
\usepackage{amsmath}

\usepackage{times}
\usepackage{bm}
\usepackage{subcaption}
\usepackage[plain]{algorithm2e}
\usepackage{graphicx}
\usepackage{amssymb}
\newcommand{\ind}{\rotatebox[origin=c]{90}{$\models$}}
\usepackage{tikz}
\usepackage{multirow}

\usepackage{soul,color}
\definecolor{mygreen}{RGB}{144,241,47}

\newtheorem{theorem}{Theorem}
\newtheorem{proposition}{Proposition}

\newtheorem{lemma}{Lemma}[]
\newtheorem{remark}{Remark}[]
\newtheorem{definition}{Definition}

\newcommand{\nind}{\not\!\perp\!\!\!\perp}
\newcommand{\mS}{\mathcal{S}}

\title{Causal inference in degenerate systems: An impossibility result}

\author{Yue Wang}
\thanks{Yue Wang, Institut des Hautes \'Etudes Scientifiques, 91440 Bures-sur-Yvette, France. Email: yuewang@ihes.fr}
\author{Linbo Wang}
\thanks{Linbo Wang, Department of Statistical Sciences, University of Toronto, Toronto, Ontario M5S 3G3, Canada. Email: linbo.wang@utoronto.ca}

\begin{document}
\begin{abstract}
Causal relationships among variables are commonly represented via directed acyclic graphs. There are many methods in the literature to quantify the strength of arrows in a causal acyclic graph. These methods, however, have undesirable properties when the causal system represented by a directed acyclic graph is degenerate. In this paper, we characterize a degenerate causal system using multiplicity of Markov boundaries. We show that in this case, it is impossible to find an identifiable quantitative measure of causal effects that satisfy a set of natural criteria. To supplement the impossibility result, we also develop algorithms to identify degenerate causal systems from observed data. Performance of our algorithms is investigated through synthetic data analysis.
\begin{flushleft}
{\bf KEY WORDS:} Causal inference; Impossibility theorem; Markov boundary.
\end{flushleft}
\end{abstract}
\maketitle

\section{Introduction}
\label{sec:introduction}

Inferring causal relationships is among the most important goals in many disciplines. A formal approach to represent causal relationships uses causal directed acyclic graphs (DAGs) (Pearl, 2009), in which random variables are represented as nodes and causal relationships are represented as arrows. Besides qualitatively describing causal relationships via DAGs, it is often desirable to obtain quantitative measures of the strength of arrows therein since they provide more detailed information on causal effects. There have been many measures proposed to quantify the causal relationships between nodes in a causal DAG, such as conditional mutual information (Dobrushin, 1963), causal strength (Janzing et al., 2013) and part mutual information (Zhao et al., 2016). See Gao et al. (2016) and its reference list for more such measures.

An interesting observation is that these measures have undesirable properties when the causal system under consideration is degenerate. As a simple example, consider the confounder triangle $Z \rightarrow X \rightarrow Y$ with an edge $Z \rightarrow Y$, where $Z=X$ almost surely. In this case, the conditional mutual information $\textsc{cmi}(X,Y\mid Z)$ is zero regardless of the influence $X$ has on $Y$, while the causal strength and part mutual information for the arrow $X\rightarrow Y$ are not well-defined. Intuitively, these problems arise because it is not possible to distinguish the causal effect of $X$ on $Y$ from the causal effect of $Z$ on $Y$. 

In this paper, we generalize the observation above by providing a formal characterization of a degenerate causal system in Section \ref{sec:multiple}. We first define a set of natural criteria to be expected from a reasonable measure of causal influence, and show that when the causal system is degenerate, all reasonable measures of a causal influence cannot be identified from the distribution represented by the DAG. Analysts may instead report qualitative summaries of causal relationships, such as all the causal explanations of the response variable. 

Our characterization of a degenerate causal system is based on multiplicity of Markov boundaries for the response variable. The Markov boundary of a variable $W$ in a variable set ${\mathcal{S}}$ is a minimal subset of ${\mathcal{S}}$, conditional on which all the remaining variables in ${\mathcal{S}}$, excluding $W$, are rendered statistically independent of $W$ (Statnikov et al., 2013). 
 {In Section \ref{sec:algorithms}, we propose novel approaches to determine the uniqueness of Markov boundary from data. Many authors have considered methods for discovery of Markov boundaries. However, the validity of their methods often requires strong assumptions (e.g. Tsamardinos \& Aliferis, 2003; Pe{\~n}a et al., 2007; Aliferis et al., 2010), some of which even imply that the response variable has a unique Markov boundary (e.g. de Morais \&
Aussem, 2010; Mani \& Cooper, 2004). Furthermore, some of these methods output all the Markov boundaries (e.g. Statnikov et al., 2013), which are not necessary for our purpose. In contrast, our novel algorithms are more robust to model assumptions and computationally more tractable. }

\section{Background}
\label{sec:background}

\subsection{Set-up}
Consider a causal DAG $\Gamma$ with vertices $\mathcal{V}$. We say $X$ is a parent of $Y$ if the path $X\rightarrow Y$ is present in $\Gamma$, and $Y$ is a descendant of $X$ if a path $X\rightarrow \cdots \rightarrow Y$ is present in $\Gamma.$ A variable is a descendant of itself, but not a parent of itself. For a variable $W$, we use $\textsc{Des}(W)$ to denote the set consisting of all descendants of $W$, and $\textsc{Pa}(W)$ to denote the set consisting of all parents of $W$.
We assume that the probability distribution $\mathfrak{p}$ over $\mathcal{V}$ is Markov with respect to $\Gamma$ in the sense that for every $W \in {\mathcal{V}}$, $W$ is independent of ${\mathcal{V}}\setminus \textsc{Des}(W)$ conditional on $\textsc{Pa}(W)$ (Spirtes et al., 2000).

We assume that we observe independent replications of $\mathcal{V}$. 
{Let $\mathcal{S}$ be all the possible parents of $Y$, namely all the variables in $\mathcal{V}$ except for $Y$ and those that are known not to be parents of $Y$. To ease presentation, in our leading case  we assume no \emph{prior} knowledge of the causal DAG so that $\mathcal{S}=\mathcal{V}\setminus \{Y\}$. Our main result also  applies to settings where one has full or partial \emph{prior} knowledge of the DAG structure. See Remark \ref{remark:difference} for more discussions. }

Let $X$ be a possible parent of $Y$ and we are interested in the causal effect of $X$ on $Y$.
Let $\mathcal{L} = \mathcal{S} \setminus \{X\}$. 
We denote the sample space of $X,Y,\mathcal{L}$ by $\mathbb{X}, \mathbb{Y}, \mathbb{L},$ respectively. 

\subsection{Measures of causal influence}
\label{sec:measures}

We now review several measures of causal influence in the literature. We only introduce their definitions in the discrete case as they are sufficient to motivate our discussions later.

\emph{Conditional mutual information.} \quad 
The conditional mutual information between $X$ and $Y$ conditional on a set $\mathcal{C}$ is defined as (Dobrushin, 1963)
\begin{flalign*}
& \quad  \textsc{cmi}(X,Y\mid \mathcal{C}) \\
&= \sum_{x \in \mathbb{X}} \sum_{y \in \mathbb{Y}} \sum_{c \in \mathbb{C}} f(x,y, c)\log\frac{f(x,y\mid c)}{f(x\mid c)f(y\mid c)},
\end{flalign*}
where $f$ is the probability density function.
It can be shown that $\textsc{cmi}(X,Y\mid \mathcal{C})=0$ if and only if $X \ind Y \mid \mathcal{C}$. When $\mathcal{C}=\emptyset$, \textsc{cmi} is known as the mutual information (\textsc{mi}): 
$$\textsc{mi}(X,Y)=\textsc{H}(X)+\textsc{H}(Y)-\textsc{H}(X,Y),$$ where $\textsc{H}$ is the Shannon entropy. { More generally, we have $\textsc{cmi}(X,Y\mid \mathcal{C})=\textsc{mi}(\{X\}\cup \mathcal{C},Y)-\textsc{mi}(Y,\mathcal{C})$.}

$\textsc{cmi}(X,Y\mid \mathcal{C})$ quantifies additional (possibly non-linear) information contained in $Y$ regarding $X$ conditional on $\mathcal{C}$.

\emph{Causal strength.} \quad 
When the full causal DAG (thus $\mathcal{L} = \textsc{Pa}(Y) \setminus \{X\}$) is known, the causal strength of $X$ on $Y$ is defined as (Janzing et al., 2013)
\begin{flalign*}
& \quad  \textsc{cs}(X\to Y) \\
&=\sum_{x\in \mathbb{X}} \sum_{y\in \mathbb{Y}} \sum_{ l\in \mathbb{L}} f(x,y, l)\log\frac{f(y\mid x, l)}{\sum\limits_{x'\in \mathbb{X}}f(y\mid x', l)f(x')}.
\end{flalign*}
\textsc{cs} is motivated to avoid the so-called underestimation problem of \textsc{cmi}: when $X$ and $Z$ are almost the same but have strong causal effect on $Y$, both $\textsc{cmi}(X,Y\mid Z)$ and $\textsc{cmi}(Z,Y\mid X)$ are very small.

\emph{Part mutual information.} \quad The part mutual information between $X$ and $Y$ conditional on $\mathcal{C}$ is defined as (Zhao et al., 2016)
\begin{flalign*}
& \quad \textsc{pmi}(X,Y\mid \mathcal{C}) \\
&=\sum_{x\in \mathbb{X}} \sum_{y\in \mathbb{Y}} \sum_{ c\in \mathbb{C}} f(x,y, c)\log\frac{f(x,y\mid c)}{f^*(x\mid c)f^*(y\mid c)},
\end{flalign*}
where $f^*(x\mid c)=\sum_{y\in \mathbb{Y}} f(x\mid y, c)f(y),
f^*(y\mid c)=\sum_{x\in \mathbb{X}} f(y\mid x, c)f(x)$. 

\textsc{pmi} solves a similar underestimation problem of \textsc{cmi}, but it is also symmetric, and definition of which does not depend on knowledge of the full DAG.

\subsection{Markov blanket and Markov boundary}
\label{sec:examples}
We now formally discuss the notion of Markov blanket and Markov boundary. 

\begin{definition}
	Suppose that ${\mathcal{T}}$ is a set of observed variables not containing $W$. 
	A subset of $\mathcal{T}$, denoted as $\mathcal{M}$, is a Markov blanket of $W$ within $\mathcal{T}$ if
	\begin{equation*}
	\label{eqn:mb}
	W \ind ({\mathcal{T}} \setminus {\mathcal{M}}) \mid \mathcal{M}.
	\end{equation*}
\end{definition}
Using the notion of mutual information, the above condition can be written as
$\textsc{cmi}(W,{\mathcal{T}}\mid \mathcal{M})=0$, or equivalently, $\textsc{mi}(W,{\mathcal{T}})=\textsc{mi}(W,\mathcal{M})$. This suggests that the Markov blanket $\mathcal{M}$ contains all the information of ${\mathcal{T}}$ on $W$.

\begin{definition}
A Markov blanket is called a Markov boundary if none of its proper subset is a Markov blanket. In other words, a Markov boundary is a minimal Markov blanket. 
\end{definition}
Markov boundary always exists. If $W \ind \mathcal{T}$, then $\emptyset$ is the Markov boundary. If no proper subset of $\mathcal{T}$ is  a Markov blanket, then $\mathcal{T}$ is the Markov boundary.

Even though Markov boundaries are minimal, in general they are not unique. For example, consider the causal DAG in Fig. \ref{ff}. Variables $X,Y,Z$ take value in $\{0,1,2\}$, while $W$ takes value in $\{0,1\}$. Both $\{X,W\}$ and $\{Z,W\}$ are Markov boundaries of $Y$ (conditioned on $\{X,W\}$, either $Z$ or $Y$ only takes one value, which implies independence), but the probability that $(X,W)=(Z,W)$ is less than one. The multiplicity of Markov boundary also implies unfaithfulness (discussed below).

\begin{figure}[!htbp]
	\centering
	\begin{tikzpicture}
	\draw  [ultra thick][->] (0.25,-1) -- (3.75,-1);
	\draw  [ultra thick][->] (0.25,-0.875) -- (1.75,-0.125);
	\draw  [ultra thick][->] (2.25,-0.125) -- (3.75,-0.875);
	\draw [ultra thick][<-] (4.25,-1) -- (5.75,-1);
	\draw (-1,-1.5) -- (7,-1.5);
	\node at (0,-1) {$Z$};
	\node at (2,0) {$X$};
	\node at (4,-1) {$Y$};
	\node at (6,-1) {$W$};
	\node at (0,-2) {0};
	\node at (2,-2) {0};
	\node at (4,-2) {0};
	\node at (6,-2.5) {0};
	\node at (0,-3) {1};
	\node at (2,-3) {1};
	\node at (4,-3) {1};
	\node at (6,-3.5) {1};
	\node at (0,-4) {2};
	\node at (2,-4) {2};
	\node at (4,-4) {2};
	\draw [thick](0.25,-2.125) -- (1.75,-2.875);
	\draw [thick](1.75,-2.125) -- (0.25,-2.875);
	\draw [thick](0.25,-2) -- (1.75,-2);
	\draw [thick](0.25,-3) -- (1.75,-3);
	\draw [thick](0.25,-4) -- (1.75,-4);
	\draw [thick](2.25,-2) -- (3.75,-2);
	\draw [thick](2.25,-2.875) -- (3.75,-2.125);
	\draw [thick](2.25,-3.875) -- (3.75,-3.125);
	\draw [thick](2.25,-4) -- (3.75,-4);
	
	\draw [thick](4.25,-1.9) -- (5.75,-2.4);
	\draw [thick](4.25,-2.1) -- (5.75,-3.4);
	\draw [thick](4.25,-2.9) -- (5.75,-2.5);
	\draw [thick](4.25,-3.1) -- (5.75,-3.5);
	\draw [thick](4.25,-3.9) -- (5.75,-2.6);
	\draw [thick](4.25,-4.1) -- (5.75,-3.6);
	\end{tikzpicture}
	\caption{A causal DAG for which variable $Y$ has multiple Markov boundaries (Statnikov et al., 2013). Combinations of values connected with lines have positive joint probabilities. For example, $\text{pr}(Z=1,X=0,Y=0,W=0) > 0$, while $\text{pr}(Z=1, X=2, Y=1, W=1) = 0$.}
	\label{ff}
\end{figure}

Unlike the confounder triangle example in Section \ref{sec:introduction}, the two Markov boundaries of $Y$ in Fig. \ref{ff} do not coincide almost surely. However, $X$ and $\{Z,W\}$ are variation dependent in the sense that there exist $x,z,w$ ($x=2,z=0,w=0$) such that $\text{pr}(X=x)>0$, $\text{pr}(Z=z,W=w)>0$, but $\text{pr}(X=x,Z=z,W=w)=0$. This variation dependence is in fact an essential property of the multiplicity of Markov boundaries.

\begin{lemma}
	\label{lemma:1}
	Let $\Theta$ denote all Markov boundaries of $Y$ in $\mathcal{T}$, where $Y \notin \mathcal{T}$. Suppose that $X \in \mathop\cup\limits_{\mathcal{M} \in \Theta} \mathcal{M} \setminus \mathop\cap\limits_{\mathcal{M} \in \Theta} \mathcal{M} $, and $\mathcal{K} =\mathcal{T}\setminus \{X\}$. Then $X$ and $\mathcal{K}$ are variation dependent in that there exist $x\in \mathbb{X}, k\in \mathbb{K}$ such that $f(x)>0$, $f(k)>0$, but $f(x,k)=0$. 
\end{lemma}

It is known in the literature (Pearl \& Paz, 1985; Pearl, 1988) {that several conditions are sufficient for the uniqueness of Markov boundary.}

\begin{definition}[Faithfulness]
	{Let $\Lambda$ denote the collection of conditional independence relationships shared by all probability distributions that are Markov with respect to $\Gamma$. A probability distribution is faithful to $\Gamma$ if and only if its conditional independence relationships are fully characterized by $\Lambda$.}
\end{definition}

\begin{definition}[Intersection]
	{A probability distribution} on ${\mathcal{V}}$ satisfies the {intersection} property if and only if for any four subsets of $\mathcal{V}$, denoted as ${\mathcal{P}}$, ${\mathcal{Q}}$, $\mathcal{Z}$, ${\mathcal{W}}$ such that ${\mathcal{P}}\ind \mathcal{Z}\mid ({\mathcal{Q}},{\mathcal{W}})$, ${\mathcal{P}}\ind {\mathcal{W}}\mid ({\mathcal{Q}},\mathcal{Z})$, it holds that ${\mathcal{P}}\ind (\mathcal{Z},{\mathcal{W}})\mid {\mathcal{Q}}$.
\end{definition}

\begin{definition}[Strict positiveness]
	{A probability distribution} on ${\mathcal{V}}$ is called {strictly positive} if and only if for any two disjoint subsets of variables ${\mathcal{X}}$ and ${\mathcal{Z}}$ such that $\text{pr}({\mathcal{X}}=x)>0$, $\text{pr}({\mathcal{Z}}={ z})>0$, it holds that $\text{pr}({\mathcal{X}}=x,{\mathcal{Z}}=z)>0$.
\end{definition}
Strict positivity allows for the expression of causal effects as conditional distributions. Nevertheless, as the proposition below shows, multiplicity of Markov boundaries implies violation of strict positivity.

\begin{proposition}[Pearl \& Paz, 1985; Pearl, 1988]
	\label{pr1}
	{If a probability distribution on $\mathcal{V}$ (i) is faithful to $\Gamma$, or (ii) has the intersection property, or (iii) is strictly positive, then any variable $Y\in \mathcal{V}$ has a unique Markov boundary in $\mathcal{V} \setminus\{Y\}$}. 
\end{proposition}

\begin{remark}
	None of the three conditions in Proposition \ref{pr1} is necessary for the uniqueness of Markov boundary. 
	{For} example, suppose $X$, $Y$, $Z$, $W$ $\in \{0,1\}$, $\text{pr}(X=Z=Y=W=0)=0.5$, $\text{pr}(X=Z=1, Y=W=0)=0.25$, $\text{pr}(X=Z=Y=W=1)=0.25$ and ${\mathcal{S}} = \{X,Y,Z,W\}$. Since $Y\ind X\mid Z$, $Y\ind Z\mid X$, $Y\not\!\perp\!\!\!\perp (X,Z)$, the joint distribution of $(X,Z,Y,W)$ does not have the intersection property and is hence not faithful (Pearl, 1988). On the other hand, $\text{pr}(X=0)>0$, $\text{pr}(Z=1)>0$ but $\text{pr}(X=0,Z=1)=0$. Hence the distribution is not strictly positive. However, each variable in this example has a unique Markov boundary within the other three variables.
\end{remark}

\subsection{Multiplicity of  Markov boundaries}
In practice, it often arises that the response variable of interest has multiple Markov boundaries. For instance, in breast cancer studies, several gene sets may have nearly the same effect for survival prediction (Ein-Dor et al., 2004), such that each of the gene sets is a Markov boundary of the survival indicator. In an extensive study, Statnikov et al. (2013) applied nine popular algorithms for learning multiple Markov boundaries to 13 benchmark data sets that cover a wide range of application domains, dimensionalities and sample sizes that are representative of practical settings. One response variable is identified for each data set. Across the nine algorithms, the frequency of reporting multiple Markov boundaries ranges from 46.2\% (6/13) to 100\%. Five out of the nine algorithms report multiple Markov boundaries in all 13 data sets. All algorithms suggest that there are multiple Markov boundaries in four out of the 13 data sets. These results suggest that a degenerate causal system (system with multiple Markov boundaries) shows up frequently in practice.

Proposition \ref{Prn} provides theoretical explanation for these empirical findings. Consider $n$ variables, each with the alphabet $\{1,...,m\}$. The joint distribution $\mathfrak{p}$ of these $n$ variables is randomly chosen from Dirichlet distribution $\text{Dir}(1,1,...,1)$.

\begin{proposition}
	\label{Prn}
For any $\epsilon,\delta>0$, when $m$ is larger than a threshold depending on $\epsilon,\delta$, and $n$ is larger than a threshold depending on $\epsilon,\delta,m$, with probability larger than $1-\delta$, we can find a probability distribution $\mathfrak{p}'$ with multiple Markov boundaries, such that the total variation distance between $\mathfrak{p}$ and $\mathfrak{p}'$ is smaller than $\epsilon$.
\end{proposition}

\begin{remark}
	Proposition \ref{Prn} concerns the measure of distributions that are at most $\epsilon$-distant from a degenerate distribution. A similar result is that the measure of $\lambda$-strong-faithful distributions is much less than one (Uhler et al. 2013). In fact, at most $\epsilon$-distant from a degenerate distribution implies $\lambda$-strong-unfaithfulness for proper $\lambda$, but not vice versa.
\end{remark}

\section{When is it possible to reasonably quantify a causal influence?}
\label{sec:multiple}

\subsection{Motivation}
\label{sec:pitfall}
We motivate our discussion in this section by generalizing our observation in the introduction. Specifically, we show that the causal effect measures introduced in Section \ref{sec:measures} may not be reasonable when the response variable $Y$ has multiple Markov boundaries within $\textsc{Pa}(Y)$. 

\begin{proposition}
	\label{prop:cmi}
	If $X \in \mathop\cup\limits_{\mathcal{M} \in \Theta} \mathcal{M} \setminus \mathop\cap\limits_{\mathcal{M} \in \Theta} \mathcal{M},$ then (i) $\textsc{cmi}(X,Y\mid \mathcal{L}) = 0$; (ii) $\textsc{cs}(X\rightarrow Y)$ and $\textsc{pmi}(X,Y\mid \mathcal{L})$ are not well-defined. Here $\Theta$ denotes all Markov boundaries of $Y$ in $\mathcal{S}$.
\end{proposition}

To solve problem (ii) in Proposition \ref{prop:cmi}, a naive solution is to assign a value in these degenerate scenarios. However, Proposition \ref{Ti} below shows that the resulting quantities cannot be continuous functions of the joint distribution of $(X,Y,\mathcal{L})$. Given a probability distribution $\mathfrak{p}'$, we use $\textsc{cs}[\mathfrak{p}'](X\to Y)$ and $\textsc{pmi}[\mathfrak{p}'](X,Y\mid \mathcal{L})$ to denote the corresponding causal strength and part mutual information.

\begin{proposition}
	\label{Ti}
	If $X\in \mathop\cup\limits_{\mathcal{M} \in \Theta} \mathcal{M}\setminus \mathop\cap\limits_{\mathcal{M} \in \Theta} \mathcal{M}$, then there exist two sequences of distributions on $(X,Y,\mathcal{L})$, denoted as $\{\mathfrak{p}_1,\mathfrak{p}_2,\ldots\}$ and $\{\mathfrak{p}'_1,\mathfrak{p}'_2,\ldots\}$, both of which converge to $\mathfrak{p}$ under the total variation distance, but $\lim_{i\to \infty} \textsc{cs}[\mathfrak{p}_i](X\to Y)\ne\lim_{i\to\infty} \textsc{cs}[\mathfrak{p}_i'](X\to Y)$. The same applies to $\textsc{pmi}(X,Y\mid \mathcal{L})$.
\end{proposition}

Proposition \ref{Ti} can be proved using Lemma \ref{lemma:1} and the following Lemma \ref{lemma:2}.

\begin{lemma}
	\label{lemma:2}
	
	Assume that there exist $x\in \mathbb{X}, l \in \mathbb{L}$ such that $f(x)>0$, $f( l)>0$, but $f(x, l)=0$. Then there exist two real numbers $g_1<g_2$, such that for any $g$ with $g_1<g<g_2$, any $\delta>0$, there exists a probability distribution $\mathfrak{p}'$ with total variation distance $\mathrm{d}(\mathfrak{p},\mathfrak{p}')<\delta$, such that $\textsc{cs}[\mathfrak{p}'](X\to Y)=g$. The same result applies to $\textsc{pmi}(X,Y\mid \mathcal{L})$.
\end{lemma}

Lemma \ref{lemma:2} is similar in flavor to the Picard's great theorem: if an analytic function $h$ has an essential singularity at a point $w$, then on any punctured neighborhood of $w$, $h(z)$ takes on all possible complex values, with at most a single exception. In this sense, \textsc{cs} and \textsc{pmi} are essentially singular at the probability distribution that implies multiple Markov boundaries for $Y$.

\subsection{Criteria for reasonable causal effect measures}

Motivated by our observations in Section \ref{sec:pitfall}, we now formally describe the criteria we expect from a reasonable measure of causal influence. We focus our discussion  on measures that  are functionals of the joint distribution of $Y$ and $\mathcal{S}$.

{C1. The strength of $X\to Y$ is a continuous function of the joint distribution of $Y$ and $\mathcal{S}$, under the total variation distance.}

C2. If there is a unique Markov boundary $\mathcal{M}$ of $Y$ within $\mathcal{S}$, and $X\notin\mathcal{M}$, then the strength of $X\to Y$ is $0$.

C3. If there is a unique Markov boundary $\mathcal{M}$ of $Y$ within $\mathcal{S}$, and $X\in\mathcal{M}$, then the absolute value of the strength of $X \to Y$ is at least $c(X,Y,\mathcal{M}\setminus\{X\})$. Here $c(X,Y,\mathcal{M}\setminus\{X\})$ is a positive constant, only depending on $X,Y,\mathcal{M}\setminus\{X\}$, such as $\textsc{cmi}(X,Y\mid \mathcal{M}\setminus\{X\})$.

We now explain why these criteria are considered natural. 

C1: Without continuity, a small perturbation on the observed distribution may lead to a big change in the effect measure. On the other hand, such a small perturbation on the observed distribution can be induced through a small perturbation on the causal system (e.g. coefficients in the structural equation models that generate the DAG). For \emph{identifiable} effect measures, a perturbation on the causal system can \emph{only} act on the effect measure through changing the observed data distribution. This suggests that a small perturbation on the underlying causal system may lead to a big change in the causal effect measure, which is undesirable.

C2: Since the unique Markov boundary contains all the information on $Y$ from $\mathcal{S}$, it is natural to say that $X$ has no causal effect on $Y$ if $X\notin \mathcal{M}$.

C3: Since any variable $X$ in the unique Markov boundary $\mathcal{M}$ of $Y$ contains non-trivial information of $Y$, it is natural to assign a positive value to the absolute value of strength of $X\rightarrow Y$. Variables outside of the unique Markov boundary should not interfere with the strength of $X\rightarrow Y$.

\subsection{An impossibility result}

We now introduce our main result in this section, which reveals the intrinsic difficulty 
to define measures of causal influence satisfying C1-C3 when multiple Markov boundaries of the response variable are present. 

Consider $\mathfrak{S}$, the set of probability distributions on $\mS \cup \{Y\}$. Choose a probability distribution $\mathfrak{p}\in \mathfrak{S}$, under which $Y$ has multiple Markov boundaries in $\mathcal{S}$, and $X$ is in at least one, but not all of such Markov boundaries. We are looking for an identifiable measure of the strength of $X\to Y$, $f:\mathfrak{S}\to\mathbb{R}$.
\begin{theorem}
	In any neighborhood $\mathfrak{N}$ of $\mathfrak{p}$ in $\mathfrak{S}$, all identifiable measures of the strength of $X\to Y$ must violate at least one of the criteria in C1 -- C3.
	\label{TH2} 
\end{theorem}

\begin{remark}
	Any two criteria among C1 to C3 are compatible with each other. For example, \textsc{cs} and \textsc{pmi} satisfy C2 and C3, a naive causal effect measure that takes a large positive constant value satisfies C1 and C3, and \textsc{cmi} satisfies C1 and C2. 
\end{remark}

To prove Theorem \ref{TH2}, we first introduce the tools that we shall use. For any random variable $X$, we define its perturbation $X^\epsilon$ to be a new random variable that coincides with $X$ with probability $1-\epsilon$, and equals an independent arbitrary noise variable $U_X$ otherwise. {For a group of variables, adding $\epsilon$-noise on one variable in the group changes the joint distribution of the whole group by at most $\epsilon$ under the total variation distance.} The following lemma shows that adding $\epsilon$-noise to $X$ will always decrease the information it has on $Y$, unless $X$ contains no information regarding $Y$.

\begin{lemma}[Strict Data Processing Inequality]
	Let ${\mathcal{S}}_1$ be a group of variables not containing $X$ or $Y$. If we add $\epsilon$-noise on $X$ to get $X^\epsilon$, then
	\begin{equation}
	\label{eqn:inequality}
	\textsc{cmi}(X^\epsilon,Y\mid {\mathcal{S}}_1)\le \textsc{cmi}(X,Y\mid {\mathcal{S}}_1),
	\end{equation}
	where the equality holds if and only if 
	\begin{equation}
	\label{eqn:equality}
	\textsc{cmi}(X,Y\mid {\mathcal{S}}_1)=0.
	\end{equation}
	\label{lemma:3}
\end{lemma}

The inequality part of Lemma \ref{lemma:3} is a special case of the data processing inequality in information theory (Cover \& Thomas, 2012). Intuitively, it states that transmitting data through a noisy channel cannot increase information, namely: garbage in, garbage out. {The original data processing inequality (Cover \& Thomas, 2012) states that the equality in \eqref{eqn:inequality} holds if and only if} 
\begin{equation}
\label{eqn:equality2}
\textsc{cmi}(X,Y\mid X^\epsilon, {\mathcal{S}}_1)=0.
\end{equation}
Condition \eqref{eqn:equality2} relies on the concrete form of noise, and thus difficult to check. In Lemma \ref{lemma:3}, {we strengthen the result by showing that \eqref{eqn:equality2} is equivalent to \eqref{eqn:equality}. This improvement is critical for the proof of }Lemma \ref{lemma: nl2}, {in which we describe how to perturb a distribution with multiple Markov boundaries for the response variable, so that in the new distribution the response variable has a unique Markov boundary}.

	\begin{lemma}
		\label{lemma: nl2}
		Assume $Y$ has multiple Markov boundaries within $\mathcal{S}$. Let $\mathcal{M}_0$ be one of them. If we add $\epsilon$-noise on each variable in $\mathcal{S} \setminus \mathcal{M}_0$, then in the new distribution, $\mathcal{M}_0$ is the unique Markov boundary. 
	\end{lemma}
	
	\textit{Proof of Theorem \ref{TH2}.}
		Assume $X$ is in Markov boundary $\mathcal{M}_1$, but not in Markov boundary $\mathcal{M}_2$. On one hand, following Lemma \ref{lemma: nl2} one may add $\epsilon$-noise on each variable in $\mathcal{S} \setminus \mathcal{M}_1$ so that $\mathcal{M}_1$ is the unique Markov boundary of $Y$. Letting $\epsilon\to 0$, criteria C1 and C3 imply that the absolute value of the strength of $X\to Y$ in the original distribution should be at least $c(X,Y, \mathcal{M}\setminus\{X\})$. On the other hand, one may also add $\epsilon$-noise on each variable in $\mathcal{S} \setminus \mathcal{M}_2$ so that $\mathcal{M}_2$ is the unique Markov boundary of $Y$. Letting $\epsilon\to 0$, criteria C1 and C2 then imply that the strength of $X\to Y$ in the original distribution should be zero. This constitutes a contradiction.

	\begin{remark}		
		\label{remark:difference}
		We note that the definition of $\mathcal{S}$ depends on knowledge of the DAG, so it is possible that one may obtain consistent estimates of a reasonable causal effect measure given the structure of the underlying DAG, but may not do so without this knowledge. For example, consider causal DAG $X_1\to X_2\to Y$ with $X_1=X_2$ almost surely. If the structure of the DAG is known \emph{a priori}, then one may define the strength of the arrow $X_2 \to Y$ by ignoring information on $X_1$. If on the other hand, one has no information on the structure of the DAG, then it is impossible to distinguish the causal effect of $X_1$ on $Y$ from the causal effect of $X_2$ on $Y$. In this case, Theorem \ref{TH2} suggests that it is impossible to obtain a reasonable quantification of the strength of the arrow $X_2 \to Y$ from data. In general, if knowledge on the DAG implies that  a variable $X$ is not a direct cause of $Y$, then one can exclude $X$ when considering the multiplicity of Markov boundaries of $Y$.
	\end{remark}
	
	\begin{remark}
	\label{remark:q}
	In the presence of multiple Markov boundaries, 
	one can report all variables that show up in at least one but not all of the Markov boundaries as ``potential causes'' of the response variable. Accuracy of such qualitative results depends on the success of algorithms that find multiple Markov boundaries. In contrast to DAG-learning, here
	one only needs  to learn the local structure around a target variable.
	\end{remark}

\section{Tests for the uniqueness of Markov boundary}
\label{sec:algorithms}

We develop a two-step procedure to test  the uniqueness of Markov boundary: (i) Find a Markov boundary for the response variable $Y$ within the observed data set $\mathcal{S}$; (ii) Decide if there exist other Markov boundaries, other than the one identified in (i). 

Methods for step (i) have been discussed extensively in the literature (Tsamardinos \& Aliferis, 2003; Pe{\~n}a et al., 2007; Aliferis et al., 2010). However, validity of existing methods typically rely on strong assumptions. For example, faithfulness is required in Aliferis et al. (2010), which implies the uniqueness of Markov boundary, and thus cannot be applied to our problem. Methods in Tsamardinos \& Aliferis (2003) and Pe{\~n}a et al. (2007) require that the joint distribution of $\mS \cup \{Y\}$ has the so-called composition property, that is, for any four subsets of $\mS \cup \{Y\}$, denoted as ${\mathcal{P}}$, ${\mathcal{Q}}$, $\mathcal{Z}$, ${\mathcal{W}}$, such that ${\mathcal{P}}\ind \mathcal{Z}\mid {\mathcal{Q}}$, ${\mathcal{P}}\ind {\mathcal{W}}\mid {\mathcal{Q}}$, it holds that ${\mathcal{P}}\ind (\mathcal{Z},{\mathcal{W}})\mid {\mathcal{Q}}$. 

To relax these assumptions, we develop Algorithm \ref{alg:2} that requires no extra assumptions on the joint distribution. Let $\Delta$ be a measure of association between two random variables, with a larger value of $\Delta$ indicating a stronger association: If two variables with $\Delta=d_1$ are dependent, then another two variables with $\Delta=d_2\ge d_1$ are also dependent. One example of $\Delta$ that we shall use in simulation studies is the conditional mutual information.

\begin{algorithm}[!htbp]
	\caption{An assumption-free algorithm for producing one Markov boundary}
	\label{alg:2}
	\vspace{-\bigskipamount}
	\begin{enumerate}
		\item \textbf{Input}\\
		\quad Joint distribution of ${\mathcal{S}}=\{X_1,\ldots,X_k\}$ and $Y$\\
		\item  {\bf Set} $\mathcal{M}_0={\mathcal{S}}$\\
		\item  {\bf Repeat}  \\
		\quad {\bf Set} $X_0=\arg\min_{X\in\mathcal{M}_0} \Delta(X,Y\mid \mathcal{M}_0\setminus\{X\})$\\
		\quad {\bf If} \ \ \ \ $X_0\ind Y\mid \mathcal{M}_0\setminus\{X_0\}$ \\
		\quad\quad  {\bf Set} $\mathcal{M}_0=\mathcal{M}_0\setminus \{X_0\}$\\
		\quad {\bf Until} $X_0\not\!\perp\!\!\!\perp Y\mid \mathcal{M}_0\setminus\{X_0\}$ 	 \\
		\item  {\bf Output} $\mathcal{M}_0$ is a Markov boundary
	\end{enumerate}
\end{algorithm}

In step 3 of Algorithm \ref{alg:2}, any tie-breaker works when there are several equal $\Delta$.

We now turn to step (ii). The key to our approach is the following necessary and sufficient condition for the uniqueness of Markov boundary.

\begin{definition}[Essential variable]
	A variable $W\in {\mathcal{S}}$ is called an {essential variable} for $Y$ if $Y \nind W \mid \mS \setminus \{W\}$. Denote the set of all essential variables by $\mathcal{E}$. 
\end{definition}
A variable $W$ is essential if it can provide additional information on $Y$, even when we have known all variables except $Y$. In Fig. \ref{ff}, $W$ is the only essential variable, since $X$ and $Z$ contain the same information on $Y$.

\begin{lemma}
	\label{lemma:nece}
	The set $\mathcal{E}$ is the intersection of all Markov boundaries of $Y$ within $\mathcal{S}$. 
\end{lemma}

\begin{theorem}
	\label{T04}
	Variable $Y$ has a unique Markov boundary within $\mS$ if and only if $\mathcal{E}$ is a Markov boundary of $Y$ within $\mS$. 
\end{theorem}

Theorem \ref{T04} provides a theoretical basis for Algorithm \ref{alg:4s} that determines if the output from Algorithm \ref{alg:2} is a unique Markov boundary. 

\begin{algorithm}[!htbp]
	\caption{A general algorithm for determining uniqueness of Markov boundary}
	\label{alg:4s}
	\vspace{-\bigskipamount}
	\begin{enumerate}
		\item \textbf{Input}\\
		\quad  Joint distribution of ${\mathcal{S}}=\{X_1,\ldots,X_k\}$ and $Y$\\
		\quad  An algorithm $\Omega$ which could produce one Markov boundary correctly\\
		\item {\bf Set} $\mathcal{M}_0=\{X_1,\ldots,X_m\}$ to be the result of Algorithm $\Omega$ on ${\mathcal{S}}$\\
		\item {\bf For} $i=1,\ldots,m$,  \\
		\quad  {\bf Set} $\mathcal{M}_i$ to be the result of Algorithm $\Omega$ on ${\mathcal{S}\setminus \{X_i\}}$
		\begin{equation*}
		\label{eqn:yind}
		\quad{\bf If}    \ \
		Y\ind \mathcal{M}_0 \mid \mathcal{M}_i  \qquad\qquad\qquad\qquad\qquad \qquad\qquad\qquad\qquad\qquad\qquad  \qquad\qquad\quad 
		\end{equation*}
		\quad\quad {\bf Output} \ \ \ \ $Y$ has multiple Markov boundaries \\
		\quad\quad  {\bf Terminate}\\
		\item {\bf Output} $Y$ has a unique Markov boundary
	\end{enumerate}
\end{algorithm}

Algorithm \ref{alg:4s} is closely related to the proposal that finds all the Markov boundaries for $Y$ in Statnikov et al. (2013). In fact, Algorithm \ref{alg:4s} can be viewed as running the proposal in Statnikov et al. (2013) until it produces two Markov boundaries or terminates. 

\begin{proposition}
\label{Pr4}
	Algorithms \ref{alg:2} and \ref{alg:4s} are sound and complete.
\end{proposition}

\begin{remark}
	\label{remark:1}
	The test in step (3) of Algorithm 2 aims to decide if $X_i$ is an essential variable. Alternatively, one may directly test
	\begin{equation}
	\label{eqn:test}
	X_i \ind Y \mid \mS \setminus \{X_i\}.
	\end{equation}
	This results in Algorithm S1 described in the supplementary material. 
	However, the conditional set $\mS \setminus \{X_i\}$ is generally very large, so that the conditional independence test for \eqref{eqn:test} may have low power.
\end{remark}

\begin{remark}
	\label{remark:2}
	A naive algorithm based on Theorem \ref{T04} involves first constructing the set of essential variables $\mathcal{E}$ in $\mS$, and then testing if $Y\ind \mS\mid \mathcal{E}$. This results in Algorithm S2 described in the supplementary material.
\end{remark}

\section{Simulation studies}
\label{sec:simulation}

We now evaluate the finite sample performance of the proposed methods. In our simulations, the response variable $Y$ and ten possible parents of $Y$, denoted as ${\mathcal{S}}=\{X_1,\ldots,X_{10}\}$, are all generated from Bernoulli distributions with mean $0.5$. We consider four settings that cover various scenarios regarding the uniqueness of Markov boundaries and the composition property of $\mathcal{S} \cup Y$. 

Setting 1: $X_1,\ldots,X_{10}$ are independent. $\text{pr}(Y=X_1)=0.8, \text{pr}(Y=X_2)=0.1$ and $\text{pr}(Y=X_3)=0.1$. In this case, $Y$ has a unique Markov boundary $\{X_1,X_2,X_3\}$.

Setting 2: Same as Setting 1, except that $X_4=X_2$. In this case, $Y$ has an additional Markov boundary $\{X_1,X_3,X_4\}$. In Settings 1 and 2, the composition property holds for $\mS \cup \{Y\}$.

Setting 3: $X_1,\ldots,X_8$ are independent. $Z=X_1+X_2 \mod 2$. $\text{pr}(Y=Z)=0.8, \text{pr}(Y=X_3)=0.1, \text{pr}(Y=X_4)=0.1, \text{pr}(X_9=X_{10}=Z)=0.95$ and $\text{pr}(X_9=X_{10}=1-Z)=0.05$. In this case, $Y$ has a unique Markov boundary: $\{X_1,X_2,X_3,X_4\}$.

Setting 4: $X_1,\ldots,X_7$ are independent. $Z=X_1+X_2 \mod 2$. $\text{pr}(Y=Z)=0.8, \text{pr}(Y=X_3)=0.1, \text{pr}(Y=X_4)=0.1, \text{pr}(X_{10}=Z)=0.95$ and $\text{pr}(X_{10}=1-Z)=0.05$. $X_8=X_1$, $X_9=X_2$. In this case, $Y$ has two Markov boundaries: $\{X_1,X_2,X_3,X_4\}$ and $\{X_3,X_4,X_8,X_9\}$. In Settings 3 and 4, the distribution of $\mS \cup \{Y\}$ violates the composition property. 

We compare the performance of the following algorithms that test the uniqueness of Markov boundaries for the response variable $Y$: (1) Alg. \ref{alg:4s}-AF: Algorithm \ref{alg:4s}, with $\Omega$ being Algorithm \ref{alg:2}; (2) Alg. \ref{alg:4s}-KI: Algorithm \ref{alg:4s}, with $\Omega$ being the KIAMB algorithm proposed in Pe{\~n}a et al. (2007), which requires the composition property; (3) Alg. S1; (4) Alg. S2. The Monte Carlo size is 500, and we report the success rates for each algorithm. In each setting, we run all four algorithms with sample size ranging from 300 to 30,000. The conditional independence test we employ is the G-test (Neapolitan, 2004) with significance level $\alpha=0.001$. All simulations are conducted with {\tt R}. Following Pe{\~n}a et al. (2007) and Statnikov et al. (2013), we choose the parameter $K$ to be $0.8$ in the KIAMB algorithm.

\begin{figure}[!h]
	\centering
	\begin{subfigure}{.38\textwidth}
		\centering
		\includegraphics[width=\linewidth]{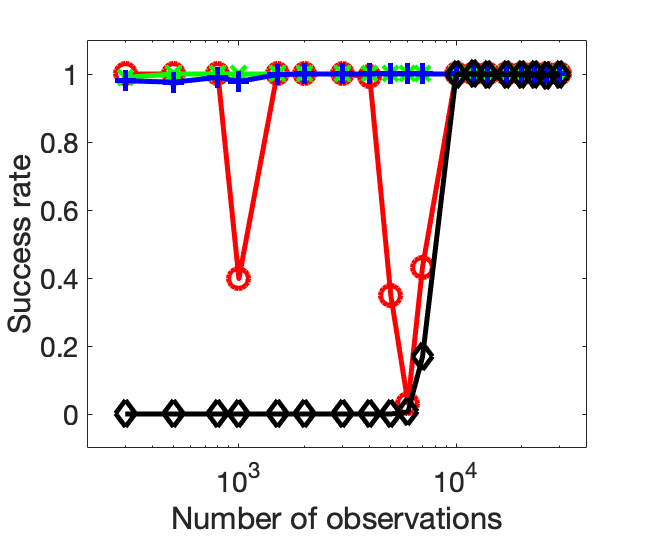}
		\caption{\small Setting 1: unique Markov boundary, composition holds.}
		\label{fig:sub1}
	\end{subfigure}
	\begin{subfigure}{.38\textwidth}
		\centering
		\includegraphics[width=\linewidth]{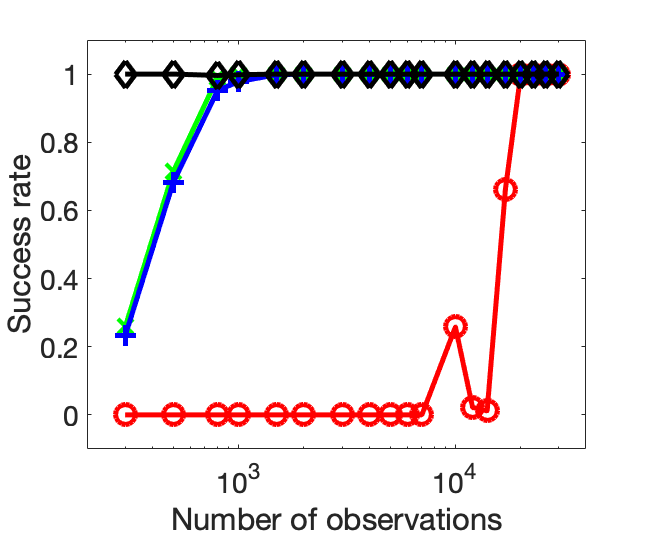}
		\caption{\small Setting 2: multiple Markov boundaries, composition holds.}
		\label{fig:sub2}
	\end{subfigure}

	\caption{Performance of various algorithms for testing the uniqueness of Markov boundary, Settings 1, 2: proposed Alg. 2-AF (blue `\textcolor{blue}{+}'); Alg. 2-KI (green `\textcolor{green}{$\times$}'); Alg. S1 (black `$\diamond$'); Alg. S2 (red `\textcolor{red}{$\circ$}'). The number of observations ranges from 300 to 30,000. The x-axis is in
logarithm scale.}
	\label{fig:simu1}
\end{figure}
\begin{figure}[!h]
	\centering

	\begin{subfigure}{.38\textwidth}
		\centering
		\includegraphics[width=\linewidth]{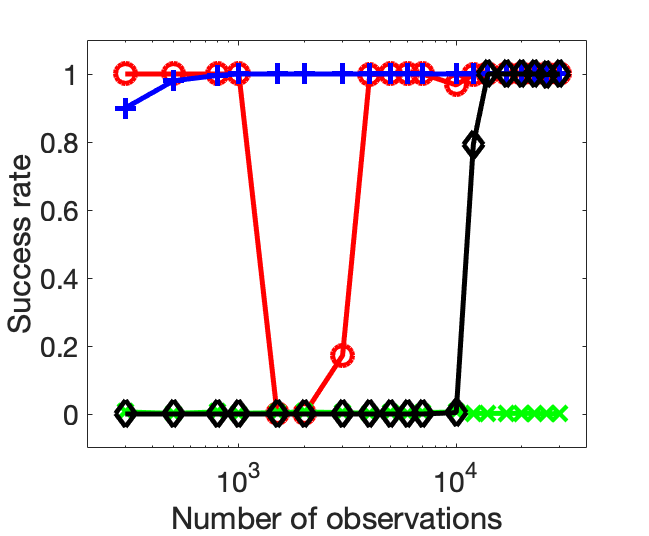}
		\caption{\small Setting 3: unique Markov boundary, composition fails.}
		\label{fig:sub3}
	\end{subfigure}
	\begin{subfigure}{.38\textwidth}
		\centering
		\includegraphics[width=\linewidth]{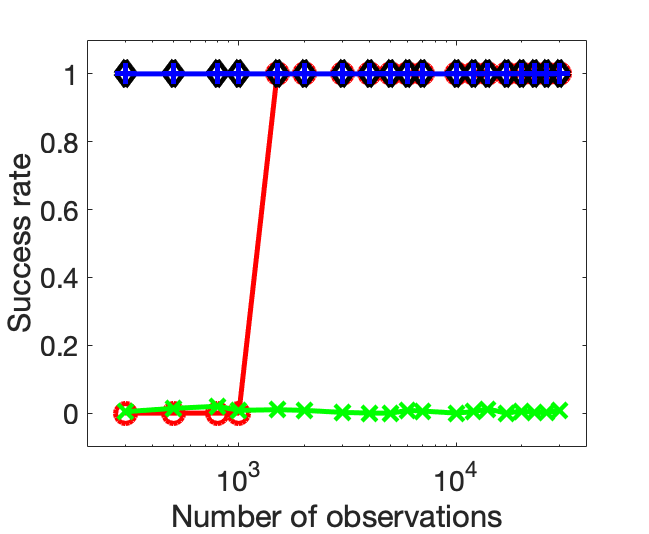}
		\caption{\small Setting 4: multiple Markov boundaries, composition fails.}
		\label{fig:sub4}
	\end{subfigure}
	\caption{Performance of various algorithms for testing the uniqueness of Markov boundary, Settings 3, 4: proposed Alg. 2-AF (blue `\textcolor{blue}{+}'); Alg. 2-KI (green `\textcolor{green}{$\times$}'); Alg. S1 (black `$\diamond$'); Alg. S2 (red `\textcolor{red}{$\circ$}'). The number of observations ranges from 300 to 30,000. The x-axis is in
logarithm scale.}
	\label{fig:simu2}
\end{figure}

As shown in Fig. \ref{fig:simu1} and Fig. \ref{fig:simu2}, both Alg. \ref{alg:4s}-AF and Alg. \ref{alg:4s}-KI have satisfactory performance under Settings 1 and 2 where the composition property holds. Alg. S1 falsely claims that there are multiple Markov boundaries for $Y$ until the sample size approaches 10,000. This is because failure to reject the hypothesis $X_i \ind Y \mid \mS \setminus \{X_i\}$ leads one to conclude that $Y$ has multiple Markov boundaries. As expected, in Settings 3 and 4 where the composition condition fails to hold, Alg. \ref{alg:4s}-AF performs much better than Alg. \ref{alg:4s}-KI. As the sample size increases, each independence test is more likely to produce correct result. When sample size is large enough, each algorithm has a high probability to produce correct final result (except Alg. \ref{alg:4s}-KI in Settings 3 and 4).

We also find that the performance of Alg. S2 is not monotonic with the number of observations. A possible explanation is that although the error rate of each single test decreases with the number of observations, certain combinations of incorrect intermediate test results  might by chance, lead to a correct final result. As the number of observations increases, the power for the independence test in step (3) of Alg. S2 increases so that the size of the empirical essential variable set $\hat{\mathcal{E}}$ grows. As a result, it is more likely that $Y\ind \mathcal{S}\mid \hat{\mathcal{E}}$ holds. On the other hand, with a larger sample size one also gains power to reject the hypothesis that $Y\ind \mathcal{S}\mid \hat{\mathcal{E}}$. This explains the non-monotonic curves we see with Alg. S2.


On average, when the composition property holds, the performance of Alg. \ref{alg:4s}-KI is slightly better than that of Alg. \ref{alg:4s}-AF, and both are much better than Alg. S1 and Alg. S2. Furthermore, Alg. \ref{alg:4s}-KI is faster than Alg. \ref{alg:4s}-AF in computation time (results not shown). When the composition property fails, Alg. \ref{alg:4s}-KI fails to produce correct results, while Alg. \ref{alg:4s}-AF exhibits the best performance.

In practice, if one has a strong belief in the composition property, then we recommend Alg. \ref{alg:4s}-KI. Otherwise Alg. \ref{alg:4s}-AF is preferable.

\subsubsection*{Acknowledgements}
The authors thank Hong Qian for motivating this paper, and Siqi He, Tengyuan Liang, Yifei Liu, Daniel Malinsky, Jifan Shi, Thomas Richardson, Weili Wang, Daxin Xu, Qingyuan Zhao and anonymous reviewers for helpful comments and discussions. Y. Wang conducted this research at the University of Washington.

	\section*{Supplementary material}

Supplementary material includes  proofs of theorems and propositions in the paper,  as well as additional algorithms referenced in Remarks \ref{remark:1} and \ref{remark:2}.

\setcounter{equation}{0}
\setcounter{figure}{0}
\setcounter{table}{0}
\makeatletter
\renewcommand{\theequation}{S\arabic{equation}}
\renewcommand{\thefigure}{S\arabic{figure}}
\renewcommand{\thetable}{S\arabic{table}}
\renewcommand{\theproposition}{S\arabic{proposition}}
\setcounter{section}{0}
\setcounter{proposition}{0}
\renewcommand{\thesection}{S}

\subsection{Proof of Lemma 1}

We first present a proposition based on the weak union property of probability distributions (Pearl, 1988).
\begin{proposition}
	Any superset of a Markov blanket is still a Markov blanket.
	\label{wup}
\end{proposition}

Now	consider two Markov boundaries $\mathcal{M}_1$, $\mathcal{M}_2$ within $\{X\}\cup \mathcal{K}$. Let $\mathcal{M}_1=\{X\}\cup \mathcal{Z}^1$, $X\notin \mathcal{M}_2$, $\mathcal{M}_2\setminus \mathcal{M}_1=\mathcal{Z}^2$, ${\mathcal{K}}\cup \{X\}\setminus (\mathcal{M}_1\cup \mathcal{M}_2)=\mathcal{Z}^3$, where $\mathcal{Z}^1=\{Z_1,\ldots,Z_n\}$, $\mathcal{Z}^2=\{ Z_1',\ldots,Z_m'\}$, $\mathcal{Z}^3=\{Z_1'',\ldots,Z_l''\}$. Therefore $\mathcal{K}=\mathcal{Z}^1\cup \mathcal{Z}^2 \cup \mathcal{Z}^3$.

Fix $z^1_0\in \mathbb{Z}^1$ such that $f (z^1_0)>0$. Assume that for $x_i\in \mathbb{X}$, $f (x_i, z^1_0)>0$ is true for $i\in\{1,\ldots,p\}$. Assume that for $ z^2_j\in \mathbb{Z}^2$, $f ( z^1_0, z^2_j)>0$ is true for $j\in\{1,\ldots,q\}$. Consider any $y\in \mathbb{Y}$.

To obtain contradiction, we assume that $f (x_i, z^1_0, z^2_j)>0$ for all $i\in\{1,\ldots,p\}$ and all $ j\in\{1,\ldots,q\}$. 

Since $X\ind Y\mid (\mathcal{Z}^1,\mathcal{Z}^2)$ (Proposition \ref{wup}) for all $i,r\in \{1,\ldots,p\}$ and all $ j\in \{1,\ldots,q\}$,
$$f (y\mid x_i, z^1_0, z^2_j)=f (y\mid x_r, z^1_0, z^2_j).$$

Since $\mathcal{Z}^2\ind Y\mid (X, \mathcal{Z}^1)$ for all $r\in \{1,\ldots,p\}$ and all $j,s\in \{1,\ldots,q\}$,
$$f (y\mid x_r, z^1_0, z^2_j)=f (y\mid x_r, z^1_0, z^2_s).$$
All the conditions have positive probabilities, so the conditional probabilities are well-defined.

Then we have 
$$f (y\mid x_i, z^1_0, z^2_j)=f (y\mid x_r, z^1_0, z^2_s),$$
for all $i,r\in \{1,\ldots,p\}$ and all $j,s\in \{1,\ldots,q\}$. 

Since this is true for any possible values of $X$ and $\mathcal{Z}^2$ when $\mathcal{Z}^1= z^1_0$, we know that 
$$f (y\mid x_i, z^1_0, z^2_j)=f (y\mid  z^1_0).$$

Therefore, for all $ z^1_1\in \mathcal{Z}^1$ with $f ( z^1_1)>0$, all $y\in \mathbb{Y}$ and all $ i,j$,
$$f (x_i, z^2_j,y\mid  z^1_1)=f (x_i, z^2_j\mid  z^1_1)f (y\mid  z^1_1)$$
is valid.

This implies that $(X, \mathcal{Z}^2)\ind Y\mid \mathcal{Z}^1$, therefore $X\ind Y\mid \mathcal{Z}^1$, $\textsc{mi}(Y,\mathcal{Z}^1)=\textsc{mi}(Y,(X,\mathcal{Z}^1))$. Since $\mathcal{M}_1=\{X\}\cup \mathcal{Z}^1$, $\textsc{mi}(Y,(X,\mathcal{Z}^1))=\textsc{mi}(Y,{\mathcal{K}})$. Thus $\textsc{mi}(Y,\mathcal{Z}^1)=\textsc{mi}(Y,\{X\}\cup {\mathcal{K}})$, implying that $\mathcal{Z}^1$ is a Markov blanket, which is a contradiction. So there exists $x\in \mathbb{X}, z^1_0\in\mathbb{Z}^1,  z^2_1\in \mathbb{Z}^2$ such that $f (x, z^1_0)>0$ (implies $f (x)>0$), $f ( z^1_0, z^2_1)>0$, but $f (x, z^1_0, z^2_1)=0$. Choose $ z^3_1\in \mathbb{Z}^3$ such that $f ( z^1_0, z^2_1, z^3_1)>0$, and let $k=( z^1_0, z^2_1, z^3_1)$, then $f (x)>0$, $f(k)>0$, but $f (x,k)=0$.

\subsection{Proof of Proposition 2}

In this setting, when $n$ is much larger than fixed $m$, due to the property of Dirichlet distribution, with probability at least $1-\delta/2$, we can modify $\mathfrak{p}$ to $\bar{\mathfrak{p}}$ such that three pre-chosen variables $X,Y,Z$ are independent under $\bar{\mathfrak{p}}$, and $\mathrm{d}(\mathfrak{p},\bar{\mathfrak{p}})<\epsilon/2$. Then construct $\bar{X},\bar{Y},\bar{Z}$: $\bar{X},\bar{Y},\bar{Z}$ equal $X,Y,Z$ if none of $X,Y,Z$ is $1$; $\bar{X},\bar{Y},\bar{Z}$ equal $1$ if at least one of $X,Y,Z$ is 1. Now either all $\bar{X},\bar{Y},\bar{Z}$ equal $1$, or none of them equals $1$ (they are independent in this case). Substitute $X,Y,Z$ by $\bar{X},\bar{Y},\bar{Z}$ to obtain a new distribution $\mathfrak{p}'$. When $m$ is large enough, $\mathrm{d}(\mathfrak{p}',\bar{\mathfrak{p}})<\epsilon/2$. Now under $\mathfrak{p}'$, $\bar{X}$ and $\bar{Z}$ contain exactly the same unique information of $\bar{Y}$, thus there exist multiple Markov boundaries. Besides, $\mathrm{d}(\mathfrak{p},\mathfrak{p}')<\epsilon/2$.

\subsection{Proof of Lemma 2}

In the following we will assume there is only one pair of $(x, l)$ such that $f (x)>0$, $f ( l)>0$, $f (x, l)=0$. If there are multiple pairs, we can treat them one by one.

We construct a family of probability distributions $\mathfrak{p}^\eta_i$ with mass functions $f^\eta_i$ based on $\mathfrak{p}$. For $(x', l')\ne (x, l)$, $f^\eta_i(x',y, l')=(1-\eta)f(x',y, l')$. $f^\eta_i(x, l)=\eta>0$, $f^\eta_i(y_j\mid x, l)=\alpha_i^j$, where $\alpha_i^j\ge 0$, $\sum_j \alpha_i^j=1$. Then for each $i$, $\textsc{cs}[\mathfrak{p}^\eta_i](X\to Y)$ can be defined, and when $\eta\to 0$, $f^\eta_i$ converges to $f$. The total variation distance between $f$ and $f_i^\eta$ is $\eta$.

When $\eta\to 0$, 
$$\textsc{cs}[\mathfrak{p}^\eta_i](X\to Y)=\sum_{x'\in \mathbb{X}} \sum_{y' \in \mathbb{Y}} \sum_{ l'\ne  l} f^\eta_i(x',y', l')\log\frac{f^\eta_i(y'\mid x',l')}{\sum_{x''\in \mathbb{X}}f^\eta_i(y'\mid x'',l')f^\eta_i(x'')} $$
$$+\sum_{x'\in \mathbb{X}} \sum_{y'\in \mathbb{Y}}  f^\eta_i(x',y', l)\log\frac{f^\eta_i(y'\mid x',l)}{\sum_{x''\in \mathbb{X}}f^\eta_i(y'\mid x'',l)f^\eta_i(x'')} $$
$$\to \sum_{x'\in \mathbb{X}} \sum_{y'\in \mathbb{Y}} \sum_{ l'\ne  l} f(x',y', l')\log\frac{f(y'\mid x',l')}{\sum_{x''\in \mathbb{X}}f(y'\mid x'',l')f(x'')}$$
$$+\sum_{x'\ne x}\sum_{y'\in \mathbb{Y}} f(x',y', l)\log f(y'\mid x', l)$$
$$-\sum_{j}f(y_j, l)\log \{f(x) \alpha_i^j +\sum_{x'\ne x}f(x') f(y_j\mid x', l)\}.$$

For different $i$, when we let $\eta\to 0$, the only different terms are 
$$-\sum_{j}f(y_j, l)\log \{ f(x) \alpha_i^j +\sum_{x'\ne x}f(x') f(y_j\mid x', l)\}.$$ We will show that the above term is not a constant with $\{\alpha_i^j\}$. Therefore we can find two groups of $\{\alpha_i^j\}$ for $i=1,2$ such that $g_1=\lim_{\eta\to 0}\textsc{cs}[\mathfrak{p}^\eta_1](X\to Y )<\lim_{\eta\to 0}\textsc{cs}[\mathfrak{p}^\eta_2](X\to Y )=g_2$.

If there is only one $y_1$ such that $f(y_1, l)>0$, then $$-\sum_{j}f(y_j, l)\log \{f(x) \alpha_i^j +\sum_{x'\ne x}f(x') f(y_j\mid x', l)\}$$ 
$$=-f(y_1, l)\log \{f(x)\alpha_i^1 +\sum_{x'\ne x}f(x') f(y_1\mid x', l)\}.$$

It is not a constant when we change $\alpha_i^1$.

If there are at least two values $y_1,y_2$ of $Y$, such that $f(y_1, l)>0$, $f(y_2, l)>0$, then we can change $\alpha_i^1$ while keeping $\alpha_i^1+\alpha_i^2=d$, and leave other $\alpha_i^j$ fixed. 

Set $f(y_1,  l)=a_1$, $f(y_2,  l)=a_2$, $f(x)=c$, $\sum_{x'\ne x}f(x') f(y_1\mid x', l)=b_1$, $\sum_{x'\ne x}f(x') f(y_2\mid x', l)=b_2$. All these terms are positive.
Then in $-\sum_{j}f(y_j, l)\log \{f(x) \alpha_i^j +\sum_{x'\ne x}f(x') f(y_j\mid x', l)\}$, terms containing $\alpha_i^1$ and $\alpha_i^2$ are 
$$-a_1\log (c \alpha_i^1 +b_1)-a_2 \log\{c(d-\alpha_i^1)+b_2\}.$$

Its derivative with respect to $\alpha_i^1$ is 
$$-\frac{a_1c}{c\alpha_i^1+b_1}+\frac{a_2c}{c(d-\alpha_i^1)+b_2}.$$

If the derivative always equal $0$ in an interval, then we should have 
$$\frac{a_1}{a_2}\equiv \frac{c\alpha_i^1+b_1}{c(d-\alpha_i^1)+b_2},$$
which is incorrect.

Now we have two groups of $\{\alpha_i^j\}$ for $i=1,2$ such that 
$$g_1=\lim_{\eta\to 0}\textsc{cs}[\mathfrak{p}^\eta_1](X\to Y )<\lim_{\eta\to 0}\textsc{cs}[\mathfrak{p}^\eta_2](X\to Y )=g_2.$$ Then for any $g\in (g_1,g_2)$, any $\delta>0$, we can find $\eta<\delta$ small enough such that $\textsc{cs}[\mathfrak{p}^\eta_1](X\to Y )<g$, $\textsc{cs}[\mathfrak{p}^\eta_2](X\to Y )>g$. Then we change $\{\alpha_1^j\}$ continuously to $\{\alpha_2^j\}$. During this process \textsc{cs} is always defined, and there exists $\{\alpha_3\}$ such that $\textsc{cs}[\mathfrak{p}^\eta_3](X\to Y)=g$.

This shows that $\textsc{cs}(X\to Y)$ is essentially ill-defined.

Since $\textsc{cs}(X\to Y)$ and $\textsc{pmi}(X,Y\mid \mathcal{L})$ have the same non-zero terms containing $f(\cdot\mid x, l)$, the same argument shows that $\textsc{pmi}(X,Y\mid \mathcal{L})$ is not well-defined.

\subsection{Proof of Lemma 3 when $X$ is discrete}

The proofs for discrete and continuous $X$ are different, therefore we state them separately. Whether $Y$ is discrete or continuous does not matter, therefore we assume $Y$ is discrete/continuous when $X$ is discrete/continuous. We impose some restrictions to simplify the proofs. If $X$ is discrete, then $U_X$ is an arbitrary discrete random variable which takes all the values of $X$ with positive probabilities. If $X$ is continuous, then $U_X$ is continuous, and its density function is always positive.

$\textsc{cmi}(X,Y\mid {\mathcal{S}}_1)=\sum_{s_1}\text{pr}(\mathcal{S}_1=s_1)\textsc{cmi}(X,Y\mid {\mathcal{S}}_1=s_1)$. For a fixed $s_1$, assume $X$  takes values $1,\ldots,r'$, $U_X$ takes values $1,\ldots,r',\ldots,r$, and $Y$ takes values $1,\ldots,t$ with positive probabilities. 
Denote $\text{pr}(X=i,Y=j\mid {\mathcal{S}}_1=s_1)$ by $p_{ij}$. Define $p_{- j}=\sum_i p_{ij}$, $p_{i -}=\sum_j p_{ij}$. With $\epsilon$-noise, $p_{- j}^\epsilon=p_{- j}$, $p_{ij}^\epsilon=(1-\epsilon)p_{ij}+\epsilon q_i p_{-j}$, $p_{i-}^\epsilon=(1-\epsilon)p_{i-}+\epsilon q_i$. Here $q_i$ is the density of $U_X$. Then we have 
$$\textsc{cmi}(X,Y\mid {\mathcal{S}}_1=s_1)=\sum_{j=1}^t \sum_{i=1}^{r'} p_{ij}\log\frac{p_{ij}}{p_{i-}p_{-j}},$$

$$\textsc{cmi}(X^\epsilon,Y\mid {\mathcal{S}}_1=s_1)=\sum_{j=1}^t \sum_{i=1}^r \{(1-\epsilon)p_{ij}+\epsilon q_i p_{-j}\}\log\frac{(1-\epsilon)p_{ij}+\epsilon q_i p_{-j}}{\{(1-\epsilon)p_{i-}+\epsilon q_i\}p_{-j}}.$$

$$\textsc{cmi}(X,Y\mid {\mathcal{S}}_1=s_1)-\textsc{cmi}(X^\epsilon,Y\mid {\mathcal{S}}_1=s_1)=$$
$$\sum_{j=1}^t \sum_{i=1}^r \Big [ (1-\epsilon+q_i \epsilon) p_{ij}\log\frac{p_{ij}}{p_{i-}p_{-j}}+\sum_{k\ne i}q_i \epsilon p_{kj}\log \frac{p_{kj}}{p_{k-}p_{-j}}$$
$$-\{(1-\epsilon)p_{ij}+\epsilon q_i p_{-j}\}\log\frac{(1-\epsilon)p_{ij}+\epsilon q_i p_{-j}}{p_{-j}\{(1-\epsilon)p_{i-}+\epsilon q_i\}} \Big ].$$
If $p_{k-}=0$, namely $k=r'+1,\ldots,r$, then we stipulate $\frac{ p_{kj}}{ p_{k-} p_{- j}}=1$. 

For fixed $i,j$ and $k=1,\ldots,r$, set 
$$a_{ij}^k=\frac{ p_{kj}}{ p_{k-} p_{- j}},$$ 
$$b_{ij}^k=\frac{\epsilon q_i p_{k-}}{(1-\epsilon)p_{i-}+\epsilon q_i} \quad \text{for} \quad k\ne i,$$
$$b_{ij}^i=\frac{(1-\epsilon+q_i \epsilon)p_{i-}}{(1-\epsilon)p_{i-}+\epsilon q_i},$$
$$c_{ij}=p_{-j}\{(1-\epsilon)p_{i-}+\epsilon q_i\}.$$
Here we know that $p_{-j}>0$, $(1-\epsilon)p_{i-}+\epsilon q_i>0$. 

Then we have
$$\textsc{cmi}(X,Y\mid {\mathcal{S}}_1=s_1)-\textsc{cmi}(X^\epsilon,Y\mid {\mathcal{S}}_1=s_1)$$
$$=\sum_{j=1}^t \sum_{i=1}^r c_{ij}\{\sum_{k=1}^r b_{ij}^ka_{ij}^k\log a_{ij}^k-(\sum_{k=1}^r a_{ij}^k b_{ij}^k)\log(\sum_{k=1}^r a_{ij}^k b_{ij}^k)\}\ge 0.$$

The last step is Jensen's inequality, since $a_{ij}^k\ge 0$, $b_{ij}^k\ge 0$, $\sum_{k=1}^r b_{ij}^k=1$, $c_{ij}> 0$, $f(x)=x\log x$ is strictly convex down when $x\ge 0$ (stipulate $0\log 0=0$). 

The equality holds if and only if for each $i,j$, $a^1_{ij}=a^2_{ij}=\cdots=a^{r'}_{ij}$, which means $p_{ij}/p_{i-}$ are equal for all $i\le r'$. Since $\sum_{i=1}^{r'} p_{i-}(p_{ij}/p_{i-})=p_{-j}$, $\sum_{i=1}^{r'} p_{i-}=1$, we have $p_{ij}/p_{i-}=p_{-j}$ for each $i,j$ such that $p_{i-}>0$ and $p_{-j}>0$. This is equivalent with that $X$ and $Y$ are independent conditioned on ${\mathcal{S}}_1=s_1$.

$\textsc{cmi}(X,Y\mid {\mathcal{S}}_1)=0$ if and only if $X$ and $Y$ are independent conditioned on any possible value of ${\mathcal{S}}_1$. Therefore, $\textsc{cmi}(X^\epsilon,Y\mid {\mathcal{S}}_1)\le \textsc{cmi}(X,Y\mid {\mathcal{S}}_1)$, and the equality holds if and only if $\textsc{cmi}(X,Y\mid {\mathcal{S}}_1)=0$.

\subsection{Proof of Lemma 3 when $X$ is continuous}

$$\textsc{cmi}(X,Y\mid {\mathcal{S}}_1)=\int_{-\infty}^\infty \textsc{cmi}(X,Y\mid {\mathcal{S}}_1=s_1) h(s_1)\mathrm{d}s_1,$$ where $h(s_1)$ is the probability density function of ${\mathcal{S}}_1$. For a fixed $s_1$, denote the joint probability density function of $X,Y$ conditioned on ${\mathcal{S}}_1=s_1$ by $p(x,y)$. Define $p_1(x)=\int_{-\infty}^\infty p(x,y)\mathrm{d}y$, $p_2(y)=\int_{-\infty}^\infty p(x,y)\mathrm{d}x$. With $\epsilon$-noise, the joint probability density function of $X,Y$ conditioned on ${\mathcal{S}}_1=s_1$ is $(1-\epsilon)p(x,y)+\epsilon q(x)p_2(y)$, where $q(x)$ is the density function of $U_X$. Notice that $\int_{-\infty}^\infty  q(x)\mathrm{d}x=1$, $\int_{-\infty}^\infty[(1-\epsilon)p(x,y)+\epsilon q(x)p_2(y)]\mathrm{d}x=p_2(y)$. Then we have 

$$\textsc{cmi}(X,Y\mid {\mathcal{S}}_1=s_1)-\textsc{cmi}(X^\epsilon,Y\mid {\mathcal{S}}_1=s_1)$$
$$=\int_{-\infty}^\infty \int_{-\infty}^\infty p(x,y)\log\frac{p(x,y)}{p_1(x)p_2(y)}\mathrm{d}x \mathrm{d}y$$
$$-\int_{-\infty}^\infty \int_{-\infty}^\infty \{(1-\epsilon)p(x,y)+\epsilon q(x)p_2(y)\}\log\frac{(1-\epsilon)p(x,y)+\epsilon q(x)p_2(y)}{\{(1-\epsilon)p_1(x)+\epsilon q(x)\}p_2(y)}\mathrm{d}x \mathrm{d}y$$
$$=\int_{-\infty}^\infty \int_{-\infty}^\infty \Big[(1-\epsilon)p(x,y)\log\frac{p(x,y)}{p_1(x)p_2(y)}$$
$$+q(x)\epsilon \Big \{\int_{-\infty}^\infty p(x_0,y)\log\frac{p(x_0,y)}{p_1(x_0)p_2(y)}\mathrm{d}x_0 \Big \}$$
$$-\{(1-\epsilon)p(x,y)+\epsilon q(x)p_2(y)\}\log\frac{(1-\epsilon)p(x,y)+\epsilon q(x)p_2(y)}{\{(1-\epsilon)p_1(x)+\epsilon q(x)\}p_2(y)} \Big ]\mathrm{d}x \mathrm{d}y.$$

For fixed $x,y$, we can define a probability measure $\mu_{x,y}(x_0)$ on $\mathbb{R}$, which is a mixture of discrete and continuous type measures. For the discrete component, it has probability $(1-\epsilon)p_1(x)/\{(1-\epsilon)p_1(x)+\epsilon q(x)\}$ to take $x$. For the continuous component, the probability density function at $x_0$ is $q(x) \epsilon p_1(x_0)/\{(1-\epsilon)p_1(x)+\epsilon q(x)\}$. Define $F_{x,y}(x_0)=p(x_0,y)/\{p_1(x_0)p_2(y)\}$. If $p_1(x_0)=0$ or $p_2(y)=0$, stipulate $F_{x,y}(x_0)=1$.

Now we have 
$$\textsc{cmi}(X,Y\mid {\mathcal{S}}_1=s_1)-\textsc{cmi}(X^\epsilon,Y\mid {\mathcal{S}}_1=s_1)$$
$$=\int_{-\infty}^\infty \int_{-\infty}^\infty \{(1-\epsilon)p_1(x)+\epsilon q(x)\}p_2(y) \Big[ \int_{-\infty}^\infty F_{x,y}(x_0)\log F_{x,y}(x_0)\mathrm{d} \mu_{x,y}(x_0)$$
$$-\Big\{ \int_{-\infty}^\infty F_{x,y}(x_0)\mathrm{d} \mu_{x,y}(x_0)\Big\}\log\Big\{ \int_{-\infty}^\infty F_{x,y}(x_0)\mathrm{d} \mu_{x,y}(x_0)\Big\}\Big]\mathrm{d}x \mathrm{d}y\ge 0.$$

The last step is the probabilistic form of Jensen's inequality, since $F_{x,y}(x_0)$ is non-negative and integrable with probability measure $\mu_{x,y}(x_0)$, $\{(1-\epsilon)p_1(x)+\epsilon q(x)\}p_2(y)> 0$ if $p_2(y)>0$, and $f(x)=x\log x$ is strictly convex down when $x\ge 0$ (stipulate $0\log 0=0$). 

The equality holds if and only if for $p_1(x_0)>0$ and $p_2(y)>0$, $F_{x,y}(x_0)$ is a constant with $x_0$, which means $p(x_0,y)/p_1(x_0)$ is a constant almost surely. Since $\int_{-\infty}^\infty p_1(x_0) p(x_0,y)/p_1(x_0)\mathrm{d}x_0=p_2(y)$, $\int_{-\infty}^\infty p_1(x_0)=1$, we have $p(x_0,y)/p_1(x_0)=p_2(y)$ for almost surely each $x_0,y$ such that $p_1(x_0)>0$ and $p_2(y)>0$. This is equivalent with that $X$ and $Y$ are independent conditioned on ${\mathcal{S}}_1=s_1$.

$\textsc{cmi}(X,Y\mid {\mathcal{S}}_1)=0$ if and only if $X$ and $Y$ are independent conditioned on any possible value of ${\mathcal{S}}_1$, except a zero-measure set. Therefore, $\textsc{cmi}(X^\epsilon,Y\mid {\mathcal{S}}_1)\le \textsc{cmi}(X,Y\mid {\mathcal{S}}_1)$, and the equality holds if and only if $\textsc{cmi}(X,Y\mid {\mathcal{S}}_1)=0$.

\subsection{Proof of Lemma 4}

Set $\mathcal{S}=\{X,Z_1,\ldots,Z_k\}$. Remember that a Markov boundary $\mathcal{M}$ is a minimal subset of $\mathcal{S}$ such that $\textsc{mi}(\mathcal{M},Y)=\textsc{mi}(\mathcal{S},Y)$. Denote ${\mathcal{S}}$ with $\epsilon$-noise on $Z_i\notin \mathcal{M}_0$ by ${\mathcal{S}}^\epsilon$. Since $\textsc{mi}(\mathcal{M}_0,Y)=\textsc{mi}({\mathcal{S}},Y)$, $\textsc{mi}(\mathcal{M}_0,Y)\le\textsc{mi}({\mathcal{S}}^\epsilon,Y)$, $\textsc{mi}({\mathcal{S}}^\epsilon,Y)\le\textsc{mi}({\mathcal{S}},Y)$, we have $\textsc{mi}({\mathcal{S}}^\epsilon,Y)=\textsc{mi}({\mathcal{S}},Y)$. Therefore, $\mathcal{M}_0$ is still a Markov boundary after adding $\epsilon$-noise. Assume in the new distribution, there is another Markov boundary, then it contains a variable with $\epsilon$-noise: $Z_i^\epsilon$. Denote this Markov boundary by $\{Z_i^\epsilon\}\cup {\mathcal{S}}_1$. Therefore, $\textsc{cmi}(Z_i^\epsilon,Y\mid {\mathcal{S}}_1)>0$. However, from Lemma 3, this implies $\textsc{cmi}(Z_i^\epsilon,Y\mid {\mathcal{S}}_1)< \textsc{cmi}(Z_i,Y\mid {\mathcal{S}}_1)$, namely $\textsc{mi}(\{Z_i^\epsilon\}\cup {\mathcal{S}}_1,Y)< \textsc{mi}(\{Z_i\}\cup {\mathcal{S}}_1,Y)$. But $\textsc{mi}(\{Z_i^\epsilon\}\cup {\mathcal{S}}_1,Y)=\textsc{mi}( {\mathcal{S}}^\epsilon,Y)=\textsc{mi}({\mathcal{S}},Y)\ge \textsc{mi}(\{Z_i\}\cup {\mathcal{S}}_1,Y)$, which is a contradiction.

\subsection{Proof of Lemma 5}

Assume there exists a Markov boundary $\mathcal{M}$ such that $W\in \mathcal{E}$, $W\notin\mathcal{M}$. Then ${\mathcal{S}}\setminus \{W\}\supset \mathcal{M}$ is a Markov blanket (Proposition \ref{wup}), and $\textsc{cmi}(Y,{\mathcal{S}}\mid {\mathcal{S}}\setminus \{W\})=0$, which contradicts to $W\in \mathcal{E}$.

If $W\notin \mathcal{E}$, then $\textsc{cmi}(Y,{\mathcal{S}}\mid {\mathcal{S}}\setminus\{W\})=0$, and ${\mathcal{S}}\setminus\{W\}$ is a Markov blanket. This Markov blanket contains a Markov boundary, which does not contain $W$.

\subsection{Proof of Theorem 2}

If Markov boundary is unique, then $\mathcal{E}$ is just the Markov boundary, therefore $\textsc{cmi}(Y,{\mathcal{S}}\mid \mathcal{E})=0$.

If $\textsc{cmi}(Y,{\mathcal{S}}\mid \mathcal{E})=0$, then $\mathcal{E}$ is a Markov blanket, which means it should contain a Markov boundary. But $\mathcal{E}$ should be contained in every Markov boundary, therefore $\mathcal{E}$ itself is a Markov boundary. $\mathcal{E}$ as a Markov boundary cannot be a proper subset of another Markov boundary, thus the only Markov boundary is $\mathcal{E}$.

\subsection{Proof of Proposition 5}

\textit{Proof that Algorithm 1 is sound and complete.}
	There exists at least one Markov boundary. The algorithm can always terminate in finite steps and produce an output. It is easy to see that the output $\mathcal{M}_0$ is a Markov blanket. In the last step of Algorithm 1, we have checked that $X_0\not\!\perp\!\!\!\perp Y\mid \mathcal{M}_0\setminus\{X_0\}$. For $X_i\in \mathcal{M}_0$, since $\Delta(X_i,Y\mid \mathcal{M}_0\setminus\{X_i\})\ge \Delta(X_0,Y\mid \mathcal{M}_0\setminus\{X_0\})$, we also have $X_i\not\!\perp\!\!\!\perp Y\mid \mathcal{M}_0\setminus\{X_i\}$. Therefore the output of Algorithm 1 is a Markov boundary.

\textit{Proof that Algorithm 2 is sound and complete.}
	The algorithm can always terminate in finite steps and produce an output. Markov boundary $\mathcal{M}_0$ is not the unique Markov boundary if and only if there exists variable $X_i\in \mathcal{M}_0$ which is not essential, namely
	\begin{equation*}
	\label{eqn:1}
	\textsc{mi}(Y,{\mathcal{S}}\setminus\{X_i\})=\textsc{mi}(Y,{\mathcal{S}}).
	\end{equation*} 
	Moreover, since $\textsc{mi}(Y,{\mathcal{S}}\setminus\{X_i\})=\textsc{mi}(Y,\mathcal{M}_i)$ and $\textsc{mi}(Y,\mathcal{M}_0)=\textsc{mi}(Y,{\mathcal{S}})$, we have $$\textsc{mi}(Y,\mathcal{M}_i)=\textsc{mi}(Y,\mathcal{M}_0),$$ or equivalently, $$\textsc{cmi} (Y,\mathcal{M}_0\mid \mathcal{M}_i)=0.$$

\subsection{Algorithms  references in Remarks 6 and 7}

We now describe Algorithms S1 and S2 that were used in the simulation studies. 

Algorithm S1 is obtained by replacing step (3) in Algorithm 2 with a direct test of  whether $X_i$ is an essential variable.

\begin{algorithm}[!htbp]
	\TitleOfAlgo{\textbf{S1.} A variant of Algorithm 2  for testing the uniqueness of Markov boundary}
	\label{alg:s2}
	\bigskip
	\begin{enumerate}
		\item \textbf{Input}\\
		\quad Joint distribution of ${\mathcal{S}}=\{X_1,\ldots,X_k\}$ and $Y$\\
		\item  {\bf Set} $\mathcal{M}_0=\{X_1,\ldots,X_m\}$ to be the result of Algorithm 1 on ${\mathcal{S}}$\\
		\item  {\bf For} $i=1,\ldots,m$,  \\
		\quad {\bf If} \ \ \ \ $X_i\ind Y\mid {\mathcal{S}}\setminus\{X_i\}$  \\
		\quad \quad {\bf Output} \ \ \ \ $Y$ has multiple Markov boundaries \\
		\quad \quad {\bf Terminate}\\
		\item  {\bf Output} $Y$ has a unique Markov boundary 	\\
	\end{enumerate}
\end{algorithm}

\textit{Proof of correctness of Algorithm S1.}
	For a Markov boundary $\mathcal{M}_0$, it is the unique Markov boundary if and only if it coincides with $\mathcal{E}$. Therefore, we only need to check whether there exists a variable $X_i\in \mathcal{M}_0$ which is not essential, namely $X_i\ind Y\mid {\mathcal{S}}\setminus\{X_i\}$.

Algorithm S2 is constructed based on Theorem 2 directly.

\begin{algorithm}[!htbp]
	\TitleOfAlgo{\textbf{S2.} A benchmark algorithm for testing  the uniqueness of Markov boundary based on Theorem 2}
	\label{alg:s1}
	\bigskip
	\begin{enumerate}
		\item  \textbf{Input}\\
		\quad Joint distribution of ${\mathcal{S}}=\{X_1,\ldots,X_k\}$ and $Y$\\
		\item {\bf Set} $\hat{\mathcal{E}}=\emptyset$\\
		\item  {\bf For} $i = 1,\ldots, k$,  \\
		\quad  {\bf If} \ \ \ \ $X_i\not\!\perp\!\!\!\perp Y\mid {\mathcal{S}}\setminus\{X_i\}$  \\
		\quad\quad $\hat{\mathcal{E}} = \hat{\mathcal{E}} \cup \{X_i\}$\\
		\item {\bf If} $Y\ind{\mathcal{S}}\mid \hat{\mathcal{E}}$ 	\\
		\quad  {\bf Output:} $Y$ has a unique Markov boundary\\
		{\bf Else} 	\\
		\quad  {\bf Output:} $Y$ has multiple Markov boundaries 
	\end{enumerate}
	
\end{algorithm}

\thispagestyle{empty}

\end{document}